\documentclass[11pt]{amsart}

\usepackage{amsfonts, amssymb, amscd}
\numberwithin{equation}{section}

\usepackage[symbol]{footmisc}

\usepackage{bm}
\usepackage{verbatim}
\usepackage{amssymb}
\usepackage{amsmath}
\usepackage{mathrsfs}
\usepackage{graphicx}
\usepackage{tikz-cd}
\usepackage{subcaption}
\usepackage{listings}
\usepackage{subfiles}
\usepackage[toc,page]{appendix}
\usepackage{mathtools}
\usepackage{comment}
\usepackage{enumerate}
\usepackage{enumitem}
\usepackage[linesnumbered,ruled]{algorithm2e}
\usepackage[all]{xy}

\usepackage{graphicx}
\graphicspath{{images/}}

\usepackage{appendix}
\usepackage{hyperref}
\newcommand{\bb}{\bm{b}}

\newcommand{\id}{\mathrm{id}}

\newcommand{\Qq}{\mathbb{Q}}

\newcommand{\Rr}{\mathbb{R}}

\newcommand{\Weil}{\operatorname{Weil}}

\newcommand{\Supp}{\operatorname{Supp}}

\newcommand{\mult}{\operatorname{mult}}

\newcommand{\Oo}{\mathcal{O}}
\newcommand{\Ii}{{\Gamma}}

\newtheorem{thm}{Theorem}[section]
\newtheorem{conj}[thm]{Conjecture}
\newtheorem{cor}[thm]{Corollary}
\newtheorem{lem}[thm]{Lemma}
\newtheorem{defn}[thm]{Definition}
\newtheorem{prop}[thm]{Proposition}
\newtheorem{ques}[thm]{Question}

\theoremstyle{definition}

\newtheorem{cons}[thm]{Construction}

\theoremstyle{definition}

\begin{document}

\title{Sarkisov program for generalized pairs}

\author{Jihao Liu}

\address{Department of Mathematics, The University of Uath, Salt Lake City, UT 84112, USA}
\email{jliu@math.utah.edu}

\begin{abstract}
In this paper we show that any two birational Mori fiber spaces of $\Qq$-factorial gklt g-pairs are connected by a finite sequence of Sarkisov links.
\end{abstract}
\date{\today}

\maketitle
\pagestyle{myheadings}\markboth{\hfill  J.Liu \hfill}{\hfill Sarkisov program for generalized pairs\hfill}

\tableofcontents

\section{Introduction}
In this paper we work over the field of complex numbers $\mathbb C$.

The minimal model program is an attempt to classify higher dimensional projective varieties. Let $X$ be a smooth complex projective variety. There are two cases: if $K_X$ is pseudo-effective, i.e. $K_X$ belongs to the closure of the big cone of divisors on $X$, then it is expected that $X$ has a minimal model $X\dashrightarrow X_{min}$. In particular, $K_{X_{min}}$ has terminal singularities, $K_{X_{min}}$ is $\mathbb Q$-Cartier and nef. The existence of such minimal models is known when $K_X$ is big, (cf. \cite{BCHM10}).

Given smooth projective variety $X$ such that $K_X$ is pseudo-effective, it is well-known that $X$ may have more than one minimal model. Nevertheless, the birational map connecting these minimal models are well-studied. For example, Kawamata (cf. \cite{Kaw08}) shows the following: for any smooth projective variety $X$ and any minimal models $X_{min}$ and $X_{min}'$ of $X$, the induced birational map $X_{min}\dashrightarrow X_{min}'$ is a composition of flops. In particular, $X_{min}$ and $X_{min}'$ are isomorphic in codimension $1$.

If $K_X$ is not pseudo-effective, by \cite{BCHM10}, we may also run a $K_X$-MMP. It is not known whether this $K_X$-MMP terminates. However, by \cite[Corollary 1.3.2]{BCHM10}, if this $K_X$-MMP is a MMP with scaling of some ample divisor, then after a sequence of flips and divisorial contractions, it terminates with a Mori fiber space $X_m\rightarrow S$. It is also well-known that the Mori fiber spaces associated to $X$ may not be unique.

\begin{ques}\label{ques: relationship mfs}
Let $X$ be a smooth projective variety such that $K_X$ is not pseudo-effective and $X_m\rightarrow S$ and $X_m'\rightarrow S'$ two Mori fiber spaces associated to $X$, i.e. there are two $K_X$-MMPs $X\dashrightarrow X_m$ and $X\dashrightarrow X_m'$ which terminates with $X_m\rightarrow S$ and $X_m'\rightarrow S'$ respectively. 

What can we say on the relationship between these two Mori fiber spaces?
\end{ques}

The \emph{Sarkisov program}, which was first introduced by Sarkisov in \cite{Sar80} and \cite{Sar82} in order to study the conic bundles of threefolds and was further developed and generalized by Corti (\cite{Cor95}), Bruno-Matsuki (\cite{BM97}) in order to study the minimal model program of klt pairs, gives us some good understanding on the question above. More precisely, the goal of the Sarkisov program is to decompose the natural birational map $X_m\dashrightarrow X_m'$ as in Question \ref{ques: relationship mfs} into four different types of birational maps, which are called the \emph{Sarkisov links}. This kind of decomposition is particularly useful when calculating the birational automorphism group of Mori fiber spaces of Fano varieties. 

Based on the minimal model program established in \cite{BCHM10}, in \cite{HM09}, Hacon and M\textsuperscript{c}Kernan have established the Sarkisov program for any klt pairs, i.e. they show the existence of such decomposition for any Sarkisov related pairs:

\begin{thm}[{\cite[Theorem 1.3]{HM09}}]\label{thm: existence sarkisov link} 
Let $(Z,\Phi)$ be a $\Qq$-factorial klt pair such that $K_X+\Phi$ is not pseudo-effective. Assume that $\phi: X\rightarrow S$ and $\psi: Y\rightarrow T$ are two Mori fiber spaces which are obtained by running two different $(K_Z+\Phi)$-MMPs. Then the induced birational map $\sigma: X\dashrightarrow Y$ is a composition of Sarkisov links.
\end{thm}

The concept of \emph{generalized pairs} was introduced in \cite{BZ16}, and it has become clear that the study of birational geometry in this category is important. For example, the most important application of generalized pairs appears in the proof of the BAB conjecture and the existence of $n$-complements (\cite{Bir16},~\cite{Bir19}).

It is a natural question to ask whether the theory of Sarkisov program holds for generalized pairs. The intuition of this question comes from the famous M$\textsuperscript{c}$Kernan-Shokurov conjecture (cf. \cite[Conjecture 1.2]{Bir16},\cite[Conjecture 1.7]{Bir18}):

\begin{conj}[M$\textsuperscript{c}$Kernan-Shokurov conjecture]\label{conj: ms conjecture}
Let $d>0$ be an integer and $\epsilon>0$ a real number. Then there is a real number $\delta>0$ depending only on $d$ and $\epsilon$ satisfying the following. Assume that $(X,B)$ is a pair and $X\rightarrow Z$ is a contraction, such that 
\begin{itemize}
    \item $(X,B)$ is $\epsilon$-lc of dimension $d$,
    \item $K_X+B\sim_{\Rr,Z}0$, and
    \item $-K_X$ is big$/Z$,
\end{itemize}
then the discriminant $\bb$-divisor $\textbf{\rm\textbf{B}}_Z$ has coefficients in $[-\infty,1-\delta]$.
\end{conj}

In the conjecture above, it was predicted that $(Z,B_Z+M_Z)$ given by the canonical bundle formula has the structure of a generalized pair. The case when $B$ is a $\Qq$-divisor is well-known, but when $B$ is an $\Rr$-divisor, this is only very recently proved by Han and Liu (cf. \cite[Corollary 1.2]{HL19}). Nevertheless, it is now natural to generalize Conjecture \ref{conj: ms conjecture} to the category of generalized pairs, which is only possible after the recent works on the generalized canonical bundle formula by Filipazzi (cf.~\cite[Theroem 1.4]{Fil18}) for generalized $\Qq$-pairs and Han-Liu (cf.~\cite[Corollary 1.2]{HL19}) for any case. Therefore, we have the following conjecture, whose $\Qq$-pair version is stated as in \cite[Conjecture 2.4]{Bir18}:

\begin{conj}\label{conj: generalized ms conjecture}
Let $d>0$ be an integer and $\epsilon>0$ a real number. Then there is a real number $\delta>0$ depending only on $d$ and $\epsilon$ satisfying the following. Assume that $f:X\rightarrow Z$ is a contraction and $(X,B+M_X)$ is a generalized pair$/Z$ with associated nef$/Z$ $\bb$-divisor $M$, such that 
\begin{itemize}
    \item $(X,B+M_X)$ is generalized $\epsilon$-lc of dimension $d$,
    \item $K_X+B+M_X\sim_{\Rr,Z}0$, and
    \item $-K_X$ is big$/Z$,
\end{itemize}
then the generalized pair given by the sub-adjunction
$$K_X+B+M_X\sim_{\Rr}f^*(K_Z+B_Z+M_Z)$$
is generalized $\delta$-lc.
\end{conj}

In this paper, we extend the theory of the Sarkisov program to the theory of generalized pairs:

\begin{thm}\label{thm: existence generalized Sarkisov link}
Assume that 
\begin{itemize}
    \item $W\rightarrow Z$ is a contraction between normal varieties,
    \item $(W,B_W+M_W)$ is a $\mathbb Q$-factorial gklt g-pair with associated nef$/Z$ $\bb$-divisor $M$, such that $K_W+B_W+M_W$ is not pseudo-effective$/Z$,
    \item $\rho_X: W\dashrightarrow X$ and $\rho_Y: W\dashrightarrow Y$ are two $(K_W+B_W+M_W)$-MMP$/Z$ such that $(\rho_X)_*(K_W+B_W+M_W)=K_X+B_X+M_X$ and $(\rho_Y)_*(K_W+B_W+M_W)=K_X+B_Y+M_Y$,
    \item $\phi_X: X\rightarrow S_X$ is a $(K_X+B_X+M_X)$-Mori fiber space$/Z$ and $\phi_Y: X\rightarrow S_Y$ is a $(K_Y+B_Y+M_Y)$-Mori fiber space$/Z$.
\end{itemize}
\begin{center}$\xymatrix{
 & W\ar@{-->}[dl]_{\rho_X}\ar@{-->}[dr]^{\rho_Y}& \\
      X \ar@{->}[d]_{\phi_X}\ar@{-->}[rr]^{f}   &  & Y\ar@{->}[d]^{\phi_Y} \\
    S_X & &S_Y }$
\end{center}
Then
\begin{enumerate}
    \item the induced birational map $f: X\dashrightarrow Y$ is given by a finite sequence of Sarkisov links$/Z$, i.e. $f$ can be written as $X_0\dashrightarrow X_1\dots\dashrightarrow X_n\cong Y$, where each $X_{i}\dashrightarrow X_{i+1}$ is a Sarkisov link$/Z$, and
    \item for any real number $\epsilon>0$, if $(W,B_W+M_W)$ is generalized $\epsilon$-lc, then $(X_i,B_i+M_{X_i})$ is generalized $\epsilon$-lc for every $i$, where each $B_i$ is the birational transform of $B_W$ on $X_i$.
\end{enumerate} 
\end{thm}

There are three additional ingredients in Theorem \ref{thm: existence generalized Sarkisov link} comparing to Theorem \ref{thm: existence sarkisov link}. Firstly, we consider \emph{gklt g-pairs} instead of klt pairs. Secondly, we indeed construct a \emph{relative Sarkisov program} since all the birational maps are over $Z$. Finally, we show that the generalized pairs constructed in the Sarkisov program do not have worse singularities comparing to $(W,B_W+M_W)$. We remark that the last two ingredients are also unknown for the klt pair case previously. Based on these three ingredients (especially the last ingredient which controls the singularities), our result may provide some positive evidence to both Conjecture \ref{conj: ms conjecture} and Conjecture \ref{conj: generalized ms conjecture}.

It is worth to mention that our proof of Theorem \ref{thm: existence generalized Sarkisov link} is quite close to the original ideas of the Sarkisov program in dimension $3$, and it is not similar to the proof of Hacon and M\textsuperscript{c}Kernan as in \cite{HM09}, where detailed results on combinatorics of ample models are applied. Indeed, the Sarkisov links constructed in our paper can be described very precisely, and can be viewed as a ``Sarkisov program with double scaling". We refer the readers to Section 4 for more details.\vspace{2mm}

\noindent\textbf{Acknowledgement.} I would like to thank Prof. Christopher D. Hacon for suggesting the question and his constant support. I would like to thank Jingjun Han for his many helpful comments, especially for his effort during the winter of 2017-2018 when he read a very early draft of this paper. The author was partially supported by NSF research grants no: DMS-1300750, DMS-1265285 and by a grant from the Simons Foundation (Award number: 256202).

\section{Notation and conventions} 
We adopt the standard notation and definitions in \cite{Sho92} and \cite{KM98}, and will freely use them.

\begin{defn}[$\bb$-divisors] Let $X$ be a normal variety. A $\bb$-$\Rr$ Cartier $\bb$-divisor ($\bb$-divisor for short) over $X$ is the choice of a projective birational morphism $Y\to X$ from a normal variety and an $\Rr$-Cartier $\mathbb R$-divisor $M$ on $Y$ up to the following equivalence: another projective birational morphism $Y'\to X$ from a normal variety and an $\Rr$-Cartier $\Rr$-divisor $M'$ defines the same $\bb$-divisor if there is a common resolution $W\to Y$ and $W\to Y'$ on which the pullback of $M$ and $M'$ coincide. If there is a choice of birational morphism $Y\rightarrow X$ such that the corresponding $\Rr$-Cartier $\Rr$-divisor $M$ is a prime divisor, the $\bb$-divisor is called \emph{prime}.

	Let $E$ be a prime $\bb$-divisor over $X$. The \emph{center} of $E$ on $X$ is the closure of its image on $X$, and is denoted by $c_X(E)$. If $c_X(E)$ is not a divisor, $E$ is called \emph{exceptional}$/X$. If $c_X(E)$ is a divisor, we say that $E$ is \emph{on} $X$. For any $\bb$-divisor $M=\sum a_iE_i$ over $X$, where $E_i$ are prime $\bb$-divisors over $X$, we define $M_X:=\sum a_ic_X(E_i)$ to be the $\Rr$-divisor where the sum is taken over all the prime $\bb$-divisors $E_i$ which are on $X$. If all the $E_i$ are on $X$, we say that $M$ is \emph{on} $X$. 
\end{defn}

\begin{defn}[Multiplicities] Let $X$ be a normal variety, $E$ a prime divisor on $X$ and $D$ an $\Rr$-divisor on $X$. We define $\mult_ED$ to be the multiplicity of $E$ along $D$. 
Let $F$ be a prime $\bb$-divisor over $X$, $B$ an $\Rr$-Cartier $\Rr$-divisor on $X$ and $\phi: Y\to X$ a birational morphism such that $F$ is on $Y$. We define $\mult_FB:=\mult_F\phi^*D$.
\end{defn}

\begin{defn}[Pairs]\label{defn: positivity}
	A \emph{pair} $(X,B)$ consists of a normal variety $X$ and an effective $\Rr$-divisor $B$ on $X$ such that $K_X+B$ is $\Rr$-Cartier.
	Let $\phi:W\to X$
	be any log resolution of $(X,B)$ and let
	$$K_W+B_W:=\phi^{*}(K_X+B).$$
	The \emph{log discrepancy} of a prime divisor $D$ on $W$ with respect to $(X,B)$ is $1-\mult_{D}B_W$ and is denoted by $a(D,X,B).$
	For any real number $\epsilon\geq 0$, we say that $(X,B)$ is \emph{lc} (resp. \emph{klt}, $\epsilon$-\emph{lc}) if $a(D,X,B)\ge0$ (resp. $>0$, $\ge\epsilon$) for every log resolution $\phi:W\to X$ as above and every prime divisor $D$ on $W$. 
	
	We say that $(X,B)$ is $\Qq$-factorial if every $\Qq$-divisor on $X$ is $\Qq$-Cartier. 
	
	For any prime $\bb$-divisor $E$ over $X$, let $Y\rightarrow X$ be a birational morphism such that $E_Y$ is a prime divisor. The \emph{log discrepancy} of $E$ with respect to $(X,B)$ is $a(E_Y,X,B)$. We say that $(X,B)$ is \emph{terminal} if $a(E,X,B)>1$ for every exceptional$/X$ prime $\bb$-divisor $E$. 
\end{defn}

\begin{defn} Let $X$ be a normal variety, $f: X\dashrightarrow Y$ a birational map, $p: W\rightarrow X$ and $q: W\rightarrow Y$ a resolution of indeterminacy of $f$, and $D$ an $\Rr$-Cartier $\Rr$-divisor on $X$ such that $D_Y:=f_*D$ is an $\Rr$-Cartier $\Rr$-divisor on $Y$. $f$ is called $D$-non-positive (resp. $D$-negative), if
\begin{itemize}
    \item $f$ does not extract any divisor, and
    \item $p^*D=q^*D_Y+E$, where $E\geq 0$ is exceptional$/Y$ (resp. $E\geq 0$ is exceptional$/Y$, and $\Supp p_*E$ contains all the $f$-exceptional divisors). 
\end{itemize}
\end{defn}

\begin{defn} We define $\Weil_{\Rr}(X)$ to be the $\Rr$-vector space spanned by all the Weil divisors on $X$. Let $\mathcal{V}$ be a finite dimensional subspace of $\Weil_{\Rr}(X)$ and $A\in\mathcal{V}$ an $\mathbb R$-divisor. We define 
$$\mathcal{L}_{A}(\mathcal{V}):=\{B|(X,B) \text{ is lc, } B=A+B', B'\geq 0, B'\in \mathcal{V}\}\subset\Weil_{\Rr}(X).$$
By \cite[Lemma 3.7.2]{BCHM10}, if $\mathcal{V}$ is a rational subspace, then $\mathcal{L}_{A}(\mathcal{V})$ is a rational polytope.
\end{defn}

\begin{defn}
A \emph{contraction} is a projective morphism $f:X\to Z$ such that $f_{*}\Oo_X=\Oo_Z$. 

For any $\bb$-divisor $M$ over $X$, $M$ is called nef$/Z$ if there is a projective morphism $Y\rightarrow X$ such that $M$ is on $Y$ and $M_Y$ is nef$/Z$.
\end{defn}

\begin{defn} A generalized pair (g-pair for short) consists of a pair $(X,B)$, a contraction $X\rightarrow Z$ and a $\bb$-divisor $M$ over $X$ such that $M$ is nef$/Z$ and $M_X$ is $\Rr$-Cartier. If there is no confusion, we usually say that $(X,B+M_X)$ is a generalized pair$/Z$. $M$ is called the associated nef$/Z$ $\bb$-divisor of the generalized pair $(X,B+M_X)$. If $Z$ is not important, we may omit $Z$ and say that $(X,B+M_X)$ is a generalized pair.

Let $(X,B+M_X)$ be a generalized pair$/Z$ with associated nef$/Z$ $\bb$-divisor $M$. Let $\phi:W\to X$
	be a log resolution of $(X,B)$ such that $M_W=M$ (i.e. $M$ is the choice of $M_W$ and the morphism $\phi$) and
	$$K_W+B_W+M_W:=\phi^{*}(K_X+B+M_X).$$
	The \emph{generalized log discrepancy} of a prime divisor $D$ on $W$ with respect to $(X,B+M_X)$ is $1-\mult_{D}B_W$ and is denoted by $ga(D,X,B+M_X).$ For any prime $\bb$-divisor $E$ over $X$, let $Y\rightarrow X$ be a birational morphism such that $E_Y$ is a prime divisor.  The \emph{generalized log discrepancy} of $E$ with respect to $(X,B+M_X)$ is $ga(E_Y,X,B+M_X)$.
	For any real number $\epsilon\geq 0$, we say that
	\begin{itemize}
	    \item $(X,B+M_X)$ is \emph{glc} (resp. \emph{gklt}, $\epsilon$-\emph{glc}) if $ga(E,X,B)\ge0$ (resp. $>0$, $\ge\epsilon$) for every prime $\bb$-divisor $E$ over $X$,
	    \item  $(X,B+M_X)$ is \emph{g-terminal} if $ga(E,X,B)>1$ for every exceptional$/X$ prime $\bb$-divisor $E$,
	    \item $(X,B+M_X)$ is \emph{gdlt} if $ga(D,X,B)>0$ for some log resolution $\phi:W\to X$ as above and every prime divisor $D$ on $W$ that is exceptional$/X$,
	    \item $(X,B+M_X)$ is $\Qq$-factorial if every $\Qq$-divisor on $X$ is $\Qq$-Cartier.
	\end{itemize}
A \textit{generalized terminalization} of a glc g-pair $(X,B+M_X)$ is a birational morphism $f: Y\rightarrow X$ satisfying the following.
\begin{itemize}
    \item $K_Y+B_Y+M_Y=f^*(K_X+B+M_X)$,
    \item $(Y,B_Y+M_Y)$ is g-terminal,
    \item $f$ only extracts prime $\bb$-divisors $E$ over $X$ such that $0\leq ga(E,X,B+M)\leq 1$.
\end{itemize}
\end{defn}

\begin{defn}
Assume that
\begin{itemize}
    \item $X\rightarrow Z$ and $Y\rightarrow Z$ are two contractions,
    \item $(X,B+M_X)$ and $(Y,B_Y+M_Y)$ be two $\Qq$-factorial generalized pairs$/Z$ with the same associated nef$/Z$ $\bb$-divisor $M$, and
    \item $f: X\dashrightarrow Y$ is a birational map$/Z$,
\end{itemize}
such that
\begin{itemize}
    \item $f$ does not extract any divisor, and
    \item $ga(E,X,B+M_X)\leq ga(E,Y,B_Y+M_Y)$ for every prime $\bb$-divisor $E$ over $X$,
\end{itemize}
then we may write $(X,B+M_X)\geq (Y,B_Y+M_Y)$.
\end{defn}

\begin{defn}[Sarkisov links]
Assume that
\begin{itemize}
    \item $X_1\rightarrow Z$ and $X_2\rightarrow Z$ are two contractions,
    \item $(X_1,B_1+M_{X_1})$ and $(X_2,B_2+M_{X_2})$ are two gklt g-pairs with the same associated nef$/Z$ $\bb$-divisor $M$,
    \item $\phi_1: X_1\rightarrow S_1$ is a $(K_{X_1}+B_1+M_{X_1})$-Mori fiber space$/Z$ and $\phi_2: X_2\rightarrow S_2$ is a $(K_{X_2}+B_2+M_{X_2})$-Mori fiber space$/Z$, 
    \item there are two birational morphisms $W\rightarrow X_1$ and $W\rightarrow X_2$ and an effective $\Rr$-divisor $B_W$ on $W$, such that $B_1$ and $B_2$ are the pushforwards of $B_W$ on $X_1$ and $X_2$ respectively, and
    \item  $f: X_1\dashrightarrow X_2$ is the induced birational map$/Z$.
\end{itemize}
Then 
\begin{itemize}
    \item $f$ is called a $(K_{X_1}+B_1+M_{X_1})$-\emph{Sarkisov link$/Z$ of type I}, or a \emph{Sarkisov link$/Z$ of type I}, if there exists an extraction $g: V\rightarrow X_1$, a sequence of flips $V\dashrightarrow X_2$ over $Z$, and an extremal contraction $S_2\rightarrow S_1$ such that the following diagram commutes:
    \begin{center}$\xymatrix{
 V\ar@{->}[d]_{g}\ar@{-->}[rr]& &X_2\ar@{->}[d]^{\phi_2} \\
      X_1\ar@{-->}[rru]_{f}\ar@{->}[dr]_{\phi_1}   &  & S_2\ar@{->}[dl] \\
    & S_1 &}$
\end{center}
\item $f$ is called a $(K_{X_1}+B_1+M_{X_1})$-\emph{Sarkisov link$/Z$ of type II}, or a \emph{Sarkisov link$/Z$ of type II}, if there exists an extraction $g: V\rightarrow X_1$, a sequence of flips $V\dashrightarrow U$ over $Z$, and a divisorial contraction $U\rightarrow X_2$, such that the following diagram commutes:
\begin{center}$\xymatrix{
 V\ar@{->}[d]_{g}\ar@{-->}[r]& U\ar@{->}[d] \\
      X_1\ar@{-->}[r]_{f}\ar@{->}[d]_{\phi_1}    & Y\ar@{->}[d]^{\phi_2} \\
    S_1\ar@{->}[r]^{\cong}& S_2 }$
\end{center}
\item $f$ is called a $(K_{X_1}+B_1+M_{X_1})$-\emph{Sarkisov link$/Z$ of type III}, or a \emph{Sarkisov link$/Z$ of type III}, if there exists a sequence of flips $X_1\dashrightarrow U$ over $Z$, a divisorial contraction $U\rightarrow X_2$ and an extremal contraction $S_1\rightarrow S_2$, such that the following diagram commutes:
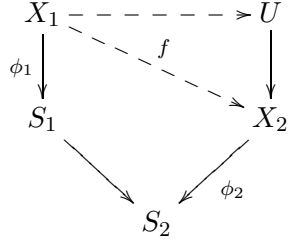
\begin{center}$\xymatrix{
 X_1\ar@{->}[d]_{\phi_1}\ar@{-->}[drr]^{f}\ar@{-->}[rr]& &U\ar@{->}[d] \\
      S_1\ar@{->}[dr]   &  & X_2\ar@{->}[dl]^{\phi_2} \\
    & S_2 &}$
\end{center}

\item $f$ is called a $(K_{X_1}+B_1+M_{X_1})$-\emph{Sarkisov link$/Z$ of type IV}, or a \emph{Sarkisov link$/Z$ of type IV}, if $f$ is a sequence of flips$/Z$, and there are two extremal contractions $S_1\rightarrow T$ and $S_2\rightarrow T$ over $Z$, such that the following diagram commutes:
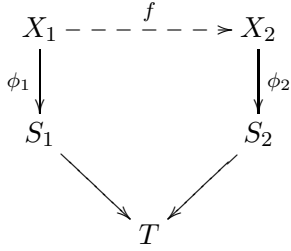
\begin{center}$\xymatrix{
 X_1\ar@{-->}[rr]^{f}\ar@{->}[d]_{\phi_1}& &X_2\ar@{->}[d]^{\phi_2} \\
      S_1\ar@{->}[dr]   &  & S_2\ar@{->}[dl] \\
    & T&}$
\end{center}
\item $f$ is called a  $(K_{X_1}+B_1+M_{X_1})$-\emph{Sarkisov link}$/Z$, or a \emph{Sarkisov link}$/Z$, if it is a Sarkisov link$/Z$ of one of the four types above. We remark that we allow $f$ to be the identity map.
\end{itemize}
\end{defn}

\section{Preliminaries} 

\begin{lem}[{\cite[Lemma 4.5 and 4.6]{BZ16}}]\label{lem: gen extraction}
Let $m\geq 0$ be an integer, $(X,B+M_X)$ a gklt g-pair, $E_1,\dots,E_m$ exceptional$/X$ $\bb$-divisors, such that $ga(E_i,X,B+M_X)\leq 1$ for each $i$. Then there is an extraction $f: Y\rightarrow X$ satisfying the following.
\begin{itemize}
\item $f$ extracts exactly $E_1,\dots,E_m$,
\item $Y$ is $\Qq$-factorial, and
\item if $m=1$ and $X$ is $\Qq$-factorial, then $f$ is a divisorial contraction.
\end{itemize}
\end{lem}

\begin{cor}\label{cor: gen terminalization}
Let $(X,B+M_X)$ be a gklt g-pair, then there is a generalized terminalization of $(X,B+M_X)$. More precisely, there is an extraction $f: Y\rightarrow X$, such that $f$ extracts exactly all the exceptional$/X$ $\bb$-divisors $E$ such that $ga(E,X,B+M_X)\leq 1$. In particular, let $K_Y+B_Y+M_Y:=f^*(K_X+B+M_X)$, then $(Y,B_Y+M_Y)$ is g-terminal.
\end{cor}

\begin{proof}
By Lemma \ref{lem: gen extraction}, we only need to show that there are finitely many prime $\bb$-divisor $E$ over $X$ such that $ga(E,X,B+M_X)\leq 1$.

Let $g: W\rightarrow X$ be a log resolution of $(X,B)$ such that $M=M_W$. Let $K_W+B_W+M_W:=g^*(K_X+B+M_X)$. Since $(X,B+M_X)$ is gklt, all the coefficients of $B_W$ are $\leq 1-c$ for some real number $c>0$. Suppose that $B_W=B_W^+-B_W^-$, where $B_W^+,B_W^-\geq 0$ and $B_W^+\wedge B_W^-=0$. Then $(W,B_W^+)$ is klt. Let $h:U\rightarrow W$ be a terminalization of $(W,B_W^+)$, then all the prime $\bb$-divisors $E$ such that $ga(E,X,B+M_X)\leq 1$ are on $U$. Therefore there are finitely many prime $\bb$-divisors $E$ over $X$ such that $ga(E,X,B+M_X)\leq 1$, and the corollary follows.
\end{proof}

\begin{prop}\label{prop: g terminalization prop} 
Let $W\rightarrow Z$ and $X\rightarrow Z$ be two contractions, $f:W\dashrightarrow X$ a birational map$/Z$, and $(W,B_W+M_W)$ and $(X,B+M_X)$ two generalized pairs$/Z$. Assume that
\begin{itemize}
    \item $K_X+B+M_X$ is nef$/Z$,
    \item $f$ does not extract any divisor,
    \item for any prime divisor $D\subset W$, $ga(D,X,B+M_X)\geq ga(D,W,B_W+M_W)$, and
    \item $(W,B_W+M_W)$ is g-terminal,
\end{itemize}
then
\begin{enumerate}
    \item $ga(E,X,B+M_X)\geq ga(E,W,B+M_W)$ for any prime $\bb$-divisor $E$ over $X$. In other words, $(W,B+M_W)\geq (X,B+M_X)$.
    \item $(X,B+M_X)$ is gklt,
    \item there is a generalized terminalization $g: Y\rightarrow X$ of $(X,B+M_X)$,
    \item the induced birational map $W\dashrightarrow Y$ does not extract any divisor,
    \item for any exceptional$/X$ $\bb$-divisor $E$ such that $ga(E,X,B+M_X)\leq 1$, $E$ is on $W$.
\end{enumerate}
\end{prop}

\begin{proof}
Let 
$p: V\rightarrow W$ and $q: V\rightarrow X$ be any resolution of indeterminacy of $f$ such that
$$p^*(K_W+B_W+M_W)=q^*(K_X+B+M_X)+E_V,$$
then $p_*E_V=\sum_{E\subset W}(ga(E,X,B+M_X)-ga(E,W,B_W+M_W))E\geq 0$. Since $K_X+B+M_X$ is nef$/Z$, $-E_V$ is nef$/W$. By the negativity lemma, $E_V\geq 0$, which implies (1). Since $(W,B_W+M_W)$ is g-terminal, $(W,B_W)$ is terminal, therefore $\lfloor B _W\rfloor=0$, and (2) follows from (1). (3) follows from (2) and Corollary \ref{cor: gen terminalization}. 

Suppose that $E$ is an exceptional$/W$ prime $\bb$-divisor. Since $(W,B_W+M_W)$ is g-terminal, $ga(E,W,B_W+M_W)>0$. Since $f$ does not extract any divisor, $E$ is exceptional$/X$. By construction of generalized terminalization, we deduce (4). (5) follows from (4).
\end{proof}

\begin{thm}[MMP for generalized pairs, {\cite[Lemma 4.4(1)(2)]{BZ16}}]\label{thm: gen pair mmp}
Let $X\rightarrow Z$ be a contraction, $(X,B+M_X)$ a $\Qq$-factorial glc g-pair$/Z$ such that $(X,0)$ is klt, and $A$ a general ample$/Z$ $\Rr$-divisor on $X$. Then there is a $(K_X+B+M_X)$-MMP$/Z$ with scaling of $A$ satisfying the following.
\begin{enumerate}
    \item If $K_X+B+M_X$ is not pseudo-effective$/Z$, then the $(K_X+B+M_X)$-MMP$/Z$ with scaling of $A$ terminates with a Mori fiber space$/Z$.
    \item If $K_X+B+M_X$ is pseudo-effective$/Z$, and
    \begin{itemize}
        \item either $(X,B+M_X)$ is gklt, or
        \item there are real numbers $\alpha,\beta\geq 0$ such that $K_X+(1+\alpha)B+(1+\beta)M_X$ is big$/Z$,
    \end{itemize}
    the $(K_X+B+M_X)$-MMP$/Z$ with scaling of $A$ terminates with a log model$/Z$ $f: X\dashrightarrow Y$, such that $K_Y+B_Y+M_Y:=f_*(K_X+B+M_X)$ is semi-ample$/Z$.
\end{enumerate}
\end{thm}

\begin{lem}\label{lem: terminalization and mmp}
Let $X\rightarrow Z$ be a contraction, $(X,B+M_X)$ a $\Qq$-factorial gklt g-pair$/Z$, and $f: X\dashrightarrow Y$ a $(K_X+B+M_X)$-non-positive map$/Z$ such that $f_*(K_X+B+M_X)=K_Y+B_Y+M_Y$. Then there is 
\begin{itemize}
    \item a resolution of indeterminacy $p: W\rightarrow X$ and $q: W\rightarrow Y$, and
    \item a g-terminal pair $(W,B_W+M_W)$,
\end{itemize}
such that
\begin{enumerate}
    \item $q$ is $(K_W+B_W+M_W)$-non-positive and $q_*(K_W+B_W+M_W)=K_Y+B_Y+M_Y$,
    \item $(W,B_W+M_W)\geq (Y,B_Y+M_Y)$,
    \item for any real number $\epsilon>0$, if $(X,B+M_X)$ is generalized $\epsilon$-lc, then $(W,B_W+M_W)$ is generalized $\epsilon$-lc, and
    \item $M=M_W$.
\end{enumerate} 
\end{lem}

\begin{proof}
Let $p: W\rightarrow X$ and $q: W\rightarrow Y$ be a resolution of indeterminacy, such that $p$ is a log resolution of $(W,B_W)$ and $M=M_W$. Let $K_W+B'_W+M_W:=p^*(K_X+B+M_X)$, then we may write
$B_W'=B_{1,W}+B_{2,W}-B_{3,W}$,
such that
\begin{itemize}
    \item $B_{1,W}$ is the strict transform of $B$ on $W$, 
    \item $B_{2,W},B_{3,W}$ are exceptional$/X$, and
    \item $B_{2,W},B_{3,W}\geq 0$ and $B_{2,W}\wedge B_{3,W}=0$.
\end{itemize}
Possibly blowing up more, we may assume that all the irreducible components of $B_{1,W}$ and $B_{2,W}$ do not intersect. Since $(X,B+M_X)$ is gklt, $(W,B_W:=B_{1,W}+B_{2,W})$ is terminal. Thus $(W,B_W+M_W)$ is g-terminal. 

Since $p$ is $(K_W+B_W+M_W)$-non-positive, $K_X+B+M_X=p_*(K_W+B_W+M_W)$ and $f$ is $(K_X+B+M_X)$-non-positive, we have that $q=f\circ p$ is $(K_W+B_W+M_W)$-non-positive, which implies (1). (2) follows from (1) and (3)(4) follow from the construction of $(W,B_W+M_W)$.
\end{proof}

\begin{thm}[Finiteness of weak log canonical models]\label{thm: finiteness ltm}
Let 
\begin{itemize}
\item $X\rightarrow Z$ be a projective morphism between normal varieties,
\item $A$ a general ample$/Z$ $\Qq$-divisor on $X$,
\item $\mathcal{V}\subset\Weil_{\Rr}(X)$ a finite dimensional rational subspace, and
\item $\mathcal{C}\subset\mathcal{L}_A(\mathcal{V})$ a rational polytope such that $(X,B)$ is klt for any $B\in\mathcal{C}$,
\end{itemize}
then there exists an integer $k\geq 0$ and birational maps$/Z$ $\phi_i: X\dashrightarrow Y_i$ for each $1\leq i\leq k$, such that for every $B\in\mathcal{C}$,
\begin{enumerate}
    \item there exists an integer $1\leq i\leq k$ such that $\phi_i$ is a weak log canonical model$/Z$ of $(X,B)$,
    \item for any log terminal model $\phi: X\dashrightarrow Y$, there exists an integer $1\leq j\leq k$ such that $\phi_j\circ\phi^{-1}: Y\rightarrow Y_j$ is an isomorphism.
\end{enumerate}
\end{thm}

\begin{proof}
The theorem follows from \cite[Theorem C]{BCHM10} and \cite[Theorem E]{BCHM10}.
\end{proof}

\section{Sarkisov program with double scaling}

\begin{thm}[Sarkisov link with double scaling]\label{thm: scaling sarkisov}
Assume that
\begin{itemize}
    \item $X\rightarrow Z$ is a contraction,
    \item $\rho: W\dashrightarrow X$ is a birational map,
    \item $(W,B_W+M_W)$ is a generalized pair with associated nef$/Z$ $\bb$-divisor $M$,
    \item $L_W$ and $H_W$ are two general big and nef$/Z$ $\Rr$-divisors on $W$, 
    \item $(X,B+M_X)$ is a generalized pair, 
    \item $\phi: X\rightarrow S$ is a $(K_X+B+M_X)$-Mori fiber space$/Z$,
    \item $\Sigma$ is a $\phi$-vertical curve,
    \item $L$ and $H$ are two $\Rr$-Cartier $\Rr$-divisors on $X$, and
    \item $0<l\leq 1$ and $0\leq h\leq 1$ are two real numbers,
\end{itemize}   
satisfying the following:
\begin{itemize}
    \item $(W,B_W+2(L_W+H_W)+M_W)$ is $\Qq$-factorial g-terminal,
    \item $K_W+B_W+H_W+M_W$ is pseudo-effective$/Z$,
    \item $(X,B+M_X)$ is $\Qq$-factorial gklt,
    \item $(W,B_W+lL_W+hH_W+M_W)\geq (X,B+lL+hH+M_X)$. In particular, $\rho$ does not extract any divisor,
    \item $B,L$ and $H$ are the birational transforms of $B_W,L_W$ and $H_W$ on $X$ respectively,
    \item $K_X+B+lL+hH+M_X\sim_{\mathbb R,S}0$, and
    \item $K_X+B+lL+hH+M_X$ is nef$/Z$.
\end{itemize}
Then there exists
\begin{itemize}
    \item a birational map$/Z$ $\rho_Y: W\dashrightarrow Y$ which does not extract any divisor,
    \item three $\Rr$-divisors $B_Y,L_Y$ and $H_Y$ on $Y$,
    \item a $(K_Y+B_Y+M_Y)$-Mori fiber space$/Z$ $\phi_Y:Y\rightarrow S_Y$,
    \item two real numbers $l\geq l_Y\geq 0$ and $h\leq h_Y\leq 1$, and
    \item a Sarkisov link$/Z$ $f: X\dashrightarrow Y$,
\end{itemize}
such that
\begin{enumerate}
\item $(Y,B_Y+M_Y)$ is a $\Qq$-factorial gklt g-pair$/Z$,
    \item $(W,B_W+l_YL_W+h_YH_W+M_W)\geq (Y,B_Y+l_YL_Y+h_YH_Y+M_Y)$. In particular, $\rho_Y$ does not extract any divisor,
    \item $B_Y,L_Y$ and $H_Y$ are the birational transforms of $B_W,L_W$ and $H_W$ on $Y$ respectively,
    \item $K_Y+B_Y+l_YL_Y+h_YH_Y+M_Y\sim_{\Rr,S}0$, 
    \item $K_Y+B_Y+l_YL_Y+h_YH_Y+M_Y$ is nef$/Z$, 
    \item for any $\phi_Y$-vertical curve $\Sigma_Y$ on $Y$,  $\frac{H_Y\cdot\Sigma_Y}{L_Y\cdot\Sigma_Y}\geq\frac{H\cdot\Sigma}{L\cdot\Sigma}>0,$ and
    \item if $h_Y=h$, then $l_Y=l$.
\end{enumerate}
\end{thm}

\begin{proof}
Since $L_W$ and $H_W$ are general big and nef$/Z$ and $\rho$ does not extract any divisor, $L$ and $H$ are big$/Z$. Thus $L\cdot\Sigma>0$ and $H\cdot\Sigma>0$, and we may define $r:=\frac{H\cdot\Sigma}{L\cdot\Sigma}$. 

For any real number $t$, we define 
$$B_W(t):=B_W+lL_W+hH_W+t(H_W-rL_W)$$ 
and 
$$B(t):=B+lL+hH+t(H-rL).$$
We define $\Ii$ to be the set of all the real numbers $t$ satisfying the following:
\begin{itemize}
    \item $0\leq t\leq\frac{l}{r}$,
        \item for any prime divisor $E\subset W$,
    $$ga(E,W,B_W(t)+M_W)\leq ga(E,X,B(t)+M_X),$$
    and
    \item $K_X+B(t)+M_X$ is nef$/Z$.
\end{itemize}
We show that $\Ii\subset [0,1-h]$. Suppose not, then there exists $t_0\in\Ii$ such that $1<h+t_0<2$. Then $(W,B_W(t_0)+M_W)$ is g-terminal, which implies that $(W,B_W(t_0)+M_W)\geq (X,B(t_0)+M_X)$. Therefore $(X,B(t_0)+M_X)$ is gklt. Since $(K_X+B(t_0)+M_X)\cdot\Sigma=0$ and $H$ is big$/Z$, 
$$(K_X+B+(l-t_0r)L+H+M_X)\cdot\Sigma=((K_X+B(t_0)+M_X)-(h+t_0-1)H)\cdot\Sigma<0.$$
Thus $\phi$ is a $(K_X+B+(l-t_0r)L+H+M_X)$-Mori fiber space$/Z$. In particular, $K_X+B+H+M_X$ is not pseudo-effective$/Z$. Since $\rho$ does not extract any divisor, $K_W+B_W+H_W+M_W$ is not pseudo-effective$/Z$, a contradiction.

Let $s:=\sup\{t|t\in\Ii\}$. In the rest of the proof, we show that we may pick $l_Y:=l-rs$ and $h_Y:=h+s$. Clearly, (7) holds under this construction. There are three cases.
\begin{itemize}
   \item[\textbf{Case 1}] $t=\frac{l}{r}$. In particular, $l_Y=0.$
    \item[\textbf{Case 2}] 
    \begin{itemize}
        \item $t<\frac{l}{r}$. In particular $l_Y>0$, and 
        \item there exists $0<\epsilon\ll 1$ and a prime divisor $E\subset W$, such that $ga(E,W,B_W(s+\epsilon)+M_W)>ga(E,X,B(s+\epsilon)+M_X).$
    \end{itemize}
        \item[\textbf{Case 3}] 
    \begin{itemize}
        \item $t<\frac{l}{r}$. In particular $l_Y>0$, and 
        \item there exists $0<\epsilon\ll 1$, such that 
        \begin{itemize}
            \item $ga(E,W,B_W(s+\epsilon)+M_W)\leq ga(E,X,B(s+\epsilon)+M_X)$ for any prime divisor $E\subset W$, and
            \item $K_X+B(s+\epsilon)+M_X$ is not nef$/Z$.
        \end{itemize} 
    \end{itemize}
\end{itemize}

\noindent\textbf{Case 1}. In this case, we finish the proof by letting $\rho_Y:=\rho, Y:=X, B_Y:=B, L_Y:=L, H_Y:=H, M_Y:=M_X, \phi_Y:=\phi_X, S_Y:=S$, and $f:=\id_X$.\vspace{2mm}

\noindent\textbf{Case 2}. In this case, $ga(E,W,B_W(s)+M_W)=ga(E,X,B(s)+M_X),$ and $E$ is exceptional$/X$. Since $E\subset W$, $$ga(E,X,B(s+\epsilon)+M_X)<ga(E,W,B_W(s+\epsilon)+M_W)\leq 1.$$
By Lemma \ref{lem: gen extraction}, there is an extraction $g: V\rightarrow X$ of $E$. By Proposition \ref{prop: g terminalization prop}, the induced birational map $W\dashrightarrow V$ does not extract any divisor. Let $0<\delta\ll\epsilon$ be a sufficiently small positive real number and
$$K_V+\Delta_V+M_V:=g^*(K_X+B+(l_Y-r\epsilon-\delta)L+(h_Y+\epsilon)H+M_X),$$
then $(V,\Delta_V+M_V)$ is a gklt g-pair. We may run a $(K_V+\Delta_V+M_V)$-MMP$/S$ $\psi: V\dashrightarrow Y$ which terminates with a Mori fiber space$/S$ $\phi_Y: Y\rightarrow S_Y$ by Theorem \ref{thm: gen pair mmp}. Since $\rho(V/S)=\rho(V/X)+\rho(X/S)=2$ and $1=\rho(Y/S_Y)\leq\rho(V/S_Y)\leq\rho(V/S)$, there are two possibilities:\vspace{2mm}

\noindent\textbf{Case 2.1}. $\rho(V/Y)=0$. In this case $\psi$ is a sequence of flips, and we get a Sarkisov link$/Z$ $f:X\dashrightarrow Y$ of type I. Let $B_Y,L_Y,H_Y$ be the birational transforms of $B_V,L_V$ and $H_V$ on $Y$ respectively and $\rho_Y: W\dashrightarrow Y$ the induced morphism. By our constructions, (1)-(5) are clear, and we only left to show (6).

For any general $\phi_Y$-vertical curve $\Sigma_Y$, $\psi$ is an isomorphism in a neighborhood of $\Sigma_Y$, and we may let $\Sigma_V$ be the birational transform of $\Sigma_Y$ on $V$. Since $\phi$ is also a $g^*(K_X+B+(l_Y-r\epsilon-\delta')L+(h_Y+\epsilon)H+M_X)$-MMP for any $0<\delta'\ll\delta$, we have
\begin{align*}
&g^*(K_X+B+(l_Y-r\epsilon-\delta')L+(h_Y+\epsilon)H+M_X)\cdot\Sigma_V\\
=&(K_Y+B_Y+(l_Y-r\epsilon-\delta')L_Y+(h_Y+\epsilon)H_Y+M_Y)\cdot\Sigma_Y<0
\end{align*}
for any $0<\delta'\leq\delta$. Thus
$$g^*(K_X+B+(l_Y-r\epsilon)L+(h_Y+\epsilon)H+M_X)\cdot\Sigma_V\leq 0.$$
Since $g^*(K_X+B+l_YL+h_YH+M_X)\sim_{\mathbb R,S}0$, we deduce that
$$g^*(H-rL)\cdot\Sigma_V\leq 0.$$
Moreover, by our assumptions, $g^*(H-rL)=g^{-1}_*(H-rL)+eE$ for some real number $e>0$, and $\Sigma_V\not\subset E$. Thus
\begin{align*}
    (H_Y-rL_Y)\cdot\Sigma_Y&=g^{-1}_*(H-rL)\cdot\Sigma_V=(g^*(H-rL)-eE)\cdot\Sigma_V\\
    &\leq g^*(H-rL)\cdot\Sigma_V\leq 0,
\end{align*}
which implies (6), and the theorem follows in this case.
\vspace{2mm}

\noindent\textbf{Case 2.2}. $\rho(V/Y)=1$. In this case, suppose that $U\rightarrow U'$ is the first divisorial contraction in the MMP$/S$: $V\dashrightarrow Y$. Then $\rho(U'/S_Y)=\rho(U'/S)=1$, which implies that $U‘\rightarrow S$ is a Mori fiber space. Thus $U'=Y$ and $S\cong S_Y$, and the induced birational map $f:X\dashrightarrow Y$ is a Sarkisov link$/Z$ of type II. Let $B_Y,L_Y,H_Y$ be the birational transforms of $B_V,L_V$ and $H_V$ on $Y$ respectively and $\rho_Y: W\dashrightarrow Y$ the induced morphism. By our constructions, (1)-(5) are clear, and we only left to show (6).

For any general $\phi_Y$-vertical curve $\Sigma_Y$, $\psi$ is an isomorphism in a neighborhood of $\Sigma_Y$, and we may let $\Sigma_V$ be the birational transform of $\Sigma_Y$ on $V$. Since $\phi$ is also a $g^*(K_X+B+(l_Y-r\epsilon-\delta')L+(h_Y+\epsilon)H+M_X)$-MMP for any $0<\delta'\ll\delta$, we have
\begin{align*}
&g^*(K_X+B+(l_Y-r\epsilon-\delta')L+(h_Y+\epsilon)H+M_X)\cdot\Sigma_V\\
=&(K_Y+B_Y+(l_Y-r\epsilon-\delta')L_Y+(h_Y+\epsilon)H_Y+M_Y)\cdot\Sigma_Y<0.
\end{align*}
for any $0<\delta'\leq\delta$. Thus
$$g^*(K_X+B+(l_Y-r\epsilon)L+(h_Y+\epsilon)H+M_X)\cdot\Sigma_V\leq 0.$$
Since $g^*(K_X+B+l_YL+h_YH+M_X)\sim_{\mathbb R,S}0$, we deduce that
$$g^*(H-rL)\cdot\Sigma_V\leq 0.$$
Moreover, by our assumptions, $g^*(H-rL)=g^{-1}_*(H-rL)+eE$ for some real number $e>0$, and $\Sigma_V\not\subset E$. Thus
\begin{align*}
    (H_Y-rL_Y)\cdot\Sigma_Y&=g^{-1}_*(H-rL)\cdot\Sigma_V=(g^*(H-rL)-eE)\cdot\Sigma_V\\
    &\leq g^*(H-rL)\cdot\Sigma_V\leq 0,
\end{align*}
which implies (6), and the theorem follows in this case.
\vspace{2mm}

\noindent\textbf{Case 3}. In this case, there exists a $(K_X+B(s+\epsilon)+M_X)$-extremal ray $[C]$ on $X$. Since $(K_X+B(s+\epsilon)+M_X)\cdot\Sigma=0$, $[C]\not=[\Sigma]$. Let $\pi: X\rightarrow T$ be the contraction of the extremal face of $\overline{NE}(X/Z)$ spanned by $[\Sigma]$ and $[C]$. Then $\pi$ factors through $S$, and  $K_X+B(s)+M_X\sim_{\Rr,T}0$.

Since $L,H$ are big$/Z$, $L,H$ are big$/T$. Therefore, if $K_{X}+B(s+\epsilon)+M_X$ is pseudo-effective$/T$, then $K_X+(1+\alpha)B(s+\epsilon)+M_X$ is big$/T$. By Theorem \ref{thm: gen pair mmp}, we may run a $(K_{X}+B(s+\epsilon)+M_X)$-MMP$/T$ with scaling of some ample $/T$ divisor, and this MMP$/T$ terminates. There are three cases:\vspace{2mm}

\noindent\textbf{Case 3.1}. After a sequence of flips $f: X\dashrightarrow Y$, the MMP$/T$ terminates with a Mori fiber space$/T$ $\phi_Y: Y\rightarrow S_Y$. Therefore, $f$ is a Sarkisov link$/Z$ of type IV.  Let $B_Y,L_Y,H_Y$ be the birational transforms of $B,L$ and $H$ on $Y$ respectively and $\rho_Y: W\dashrightarrow Y$ the induced morphism. By our constructions, (1)-(5) are clear, and we only left to show (6).

For any general $\phi_Y$-vertical curve $\Sigma_Y$, $f$ is an isomorphism in a neighborhood of $\Sigma_Y$, and we may let $\Sigma_X$ be the birational transform of $\Sigma_Y$ on $X$. Since $\phi_Y$ is a $(K_Y+B_Y+(l_Y-r\epsilon)L_Y+(h_Y+\epsilon)H_Y+M_Y)$-Mori fiber space$/T$, 
$$-(K_Y+B_Y+(l_Y-r\epsilon)L_Y+(h_Y+\epsilon)H_Y+M_Y)\cdot\Sigma_Y>0,$$
which implies that
$$-(K_X+B(s+\epsilon)+M_X)\cdot\Sigma_X>0.$$
Since $K_X+B(s)+M_X\sim_{\Rr,T}0$, 
$$-(K_X+B(s)+M_X)\cdot\Sigma_X=0,$$
which implies that 
$$(H_Y-rL_Y)\cdot\Sigma_Y=(H-rL)\cdot\Sigma_X<0.$$ 
Thus $r>\frac{H_Y\cdot\Sigma_Y}{L_Y\cdot\Sigma_Y}$, which implies (6), and the theorem follows in this case.
\vspace{2mm}

\noindent\textbf{Case 3.2}. After a sequence of flips $X\dashrightarrow U$, we get a divisorial contraction$/T$: $U\rightarrow Y$. Therefore $\rho(Y/T)=1$, which implies that the induced morphism $\phi_Y:=Y\rightarrow T$ is a Mori fiber space, and the induced birational map $f: X\dashrightarrow Y$ is a Sarkisov link$/Z$ of type III. Let $B_Y,L_Y,H_Y$ be the birational transforms of $B,L$ and $H$ on $Y$ respectively and $\rho_Y: W\dashrightarrow Y$ the induced morphism. By our constructions, (1)-(5) are clear, and we only left to show (6).

For any general $\phi_Y$-vertical curve $\Sigma_Y$, $f$ is an isomorphism in a neighborhood of $\Sigma_Y$, and we may let $\Sigma_X$ be the birational transform of $\Sigma_Y$ on $X$. Since $-(K_X+B(s+\epsilon)+M_X)$ is nef$/T$ and $K_X+B(s)+M_X\sim_{\Rr,T}0$, we have
$$-(K_X+B(s+\epsilon)+M_X)\cdot\Sigma_X\geq 0=-(K_X+B(s)+M_X)\cdot\Sigma_X,$$
which implies that 
$$(H_Y-rL_Y)\cdot\Sigma_Y=(H-rL)\cdot\Sigma_X\leq 0.$$ 
Thus $r\geq\frac{H_Y\cdot\Sigma_Y}{L_Y\cdot\Sigma_Y}$, which implies (6), and the theorem follows in this case.
\vspace{2mm}

\noindent\textbf{Case 3.3}. After a sequence of flips $f: X\dashrightarrow Y$, the MMP terminates with a minimal model $Y$ over $T$. Let $B_Y,L_Y,H_Y$ be the birational transforms of $B,L$ and $H$ on $Y$ respectively. Since $\Sigma$ is a general $\phi$-vertical curve, we may let $\Sigma'$ be the birational transform of $\Sigma$ on $Y$. Since $(K_{X}+B(s+\epsilon)+M_X)\cdot\Sigma=0$ and $(K_{X}+B(s)+M_X)\cdot\Sigma=0$, we have
$$(K_Y+B_Y+(l_Y-r\epsilon)L_Y+(h_Y+\epsilon)H_Y+M_Y)\cdot\Sigma'=0$$
and
$$(K_Y+B_Y+l_YL_Y+h_YH_Y+M_Y)\cdot\Sigma'=0$$
which implies that $(K_Y+B_Y+M_Y)\cdot\Sigma'<0$ and $r=\frac{H_Y\cdot\Sigma'}{L_Y\cdot\Sigma'}$. Let $\phi_Y:Y\rightarrow S_Y$ be the contraction of $[\Sigma']$, then $\phi_Y$ is a $(K_Y+B_Y+M_Y)$-Mori fiber space$/T$. Thus $f$ is a Sarkisov link$/Z$ of type IV.  We finish the proof by letting $\rho_Y: W\dashrightarrow Y$ be the induced birational map.
\end{proof}

\begin{lem}\label{lem: sar scaling terminalization}
Assumptions and notation as in Theorem \ref{thm: scaling sarkisov}.
\begin{enumerate}
    \item In \textbf{\rm\textbf{Case 2.1}}, $\rho(Y)-\rho(X)=1$.
    \item In \textbf{\rm\textbf{Case 2.2}},
    \begin{enumerate}
        \item $\rho(X)=\rho(Y)$,
        \item there is a prime divisor $F_0$ on $W$, such that 
    \begin{align*}
        &ga(F_0,X,B+(l_Y-r\epsilon-\delta)L+(h_Y+\epsilon)H+M_X)\\
        <&ga(F_0,Y,B_Y+(l_Y-r\epsilon-\delta)L_Y+(h_Y+\epsilon)H_Y+M_Y),
    \end{align*}
    and
    \item for any prime divisor $F$ on $W$,
      \begin{align*}
        &ga(F,X,B+(l_Y-r\epsilon-\delta)L+(h_Y+\epsilon)H+M_X)\\
        \leq &ga(F,Y,B_Y+(l_Y-r\epsilon-\delta)L_Y+(h_Y+\epsilon)H_Y+M_Y).
    \end{align*}
    \end{enumerate}
    \item In \textbf{\rm\textbf{Case 3}}, $$ga(F,Y,B_Y+(l_Y-r\epsilon)L_Y+(h_Y+\epsilon)H_Y+M_Y)\geq ga(F,W,B_W(s+\epsilon)+M_W).$$
    \item In \textbf{\rm\textbf{Case 3.1}}, $\frac{H\cdot\Sigma}{L\cdot\Sigma}>\frac{H_Y\cdot\Sigma_Y}{L_Y\cdot\Sigma_Y}$.
    \item In \textbf{\rm\textbf{Case 3.2}}, $\rho(X)-\rho(Y)=1$.
    \item In \textbf{\rm\textbf{Case 3.3}},
    \begin{enumerate}
        \item $\rho(X)=\rho(Y)$,
        \item there is a prime divisor $F_0$ on $W$, such that 
   $$ga(F_0,Y,B_Y+(l_Y-r\epsilon)L_Y+(h_Y+\epsilon)H_Y+M_Y)> ga(F_0,X,B(s+\epsilon)+M_X),$$
    and
    \item for any prime divisor $F$ on $W$,
   $$ga(F,Y,B_Y+(l_Y-r\epsilon)L_Y+(h_Y+\epsilon)H_Y+M_Y)\geq ga(F,X,B(s+\epsilon)+M_X).$$
    \end{enumerate}
\end{enumerate}
\end{lem}
\begin{proof}
(1)(4)(5) immediately follow from the proof of Theorem \ref{thm: scaling sarkisov}. (2) follows from the fact that in \textbf{\rm\textbf{Case 2.2}}, the Sarkisov link$/Z$ is constructed by running a $g^*(K_X+B+(l_Y-r\epsilon-\delta)L+(h_Y+\epsilon)H+M_X)$-MMP$/S$ and $X\not\cong Y$. (3)(6) follow from the fact that in \textbf{\rm\textbf{Case 3}}, the Sarkisov link$/Z$ is constructed by running a $(K_X+B(s+\epsilon)+M_X)$-MMP$/T$ and $X\not\cong Y$ in \textbf{\rm\textbf{Case 3.3}}.
\end{proof}

\begin{cons}[Sarkisov program with double scaling]\label{cons: sarkisov scaling}
Assume that $W\rightarrow Z$ is a contraction and $(W,B_W+M_W)$ is a $\Qq$-factorial gklt g-pair$/Z$ with nef$/Z$ $\bb$-divisor $M=M_W$, such that $K_W+B_W+M_W$ is not pseudo-effective$/Z$.

By Theorem \ref{thm: gen pair mmp}, we may run a $(K_W+B_W+M_W)$-MMP$/Z$ $\rho: W\dashrightarrow X$ which terminates with a Mori fiber space $\phi: X\rightarrow S$. 

Let $B$ be the birational transform of $B_W$ on $X$. Then $\phi$ is a $(K_{X}+B+M_{X})$-Mori fiber space$/Z$. In particular, $-(K_{X}+B+M_{X})$ is ample$/S$. Therefore we may pick a general ample$/Z$ $\Rr$-divisor $A$ on $S$ such that $-(K_{X}+B+M_{X})+\phi^*A$ is ample$/Z$. We let $L$ be a general element of $|-(K_{X}+B+M_{X})+\phi^*A|_{\Rr/Z}$ and $L_W:=\rho^*L=(\rho^{-1})_*L$. Then $L_W$ is big and nef$/Z$, $K_{X}+B+L+M_X\sim_{\Rr,S}0$ and $K_{X}+B+L+M_{X}$ is nef$/Z$.

For any general big and nef$/Z$ $\Rr$-divisor $H_W$ on $W$ such that
\begin{itemize}
\item $(W,B_W+2(L_W+H_W)+M_W)$ is $\Qq$-factorial g-terminal, and
    \item $K_W+B_W+H_W+M_W$ is pseudo-effective$/Z$, 
\end{itemize}
we construct the \emph{Sarkisov program$/Z$ of $(X,B+M_{X})$ with scaling of $(L_W,H_W)$} in the following way. 
\begin{itemize}
    \item[\textbf{Step 1}] We define $X_0:=X,B_0:=B,\rho_0:=\rho,\phi_0:=\phi$, $L_0:=L, H_0:=\rho_*H$, and $(l_0,h_0):=(1,0)$.
    \item[\textbf{Step 2}] For any integer $i\geq 0$, suppose that we already have
    \begin{itemize}
        \item a $\Qq$-factorial gklt g-pair $(X_i,B_i+M_{X_i})$,
        \item a birational map $\rho_i: W\dashrightarrow X_i$,
        \item a $(K_{X_i}+B_i+M_{X_i})$-Mori fiber space$/Z$ $\phi_i: X_i\rightarrow S_i$,
        \item two $\Rr$-Cartier $\Rr$-divisors $L_i$ and $H_i$ on $X_i$, and
            \item two real number $0<l_i\leq 1$ and $0\leq h_i\leq 1$,
    \end{itemize}
    such that
    \begin{itemize}
        \item $(W,B_W+l_iL_W+h_iH_W+M_W)\geq (X_i,B_i+l_iL_i+h_iH_i+M_{X_i})$,
        \item $B_i,L_i$ and $H_i$ are the birational transforms of $B_i,L_i$ and $H_i$ on $X_i$ respectively,
        \item $K_{X_i}+B_i+l_iL_i+h_iH_i+M_{X_i}\sim_{\Rr,S_i}0$, and
        \item $K_{X_i}+B_i+l_iL_i+h_iH_i+M_{X_i}$ is nef$/Z$,
    \end{itemize}
    then by Theorem \ref{thm: scaling sarkisov}, there exists 
    \begin{itemize}
        \item a $\Qq$-factorial gklt g-pair $(X_{i+1},B_{i+1}+M_{X_{i+1}})$,
        \item a birational map $\rho_{i+1}: W\dashrightarrow X_{i+1}$,
        \item a $(K_{X_{i+1}}+B_{i+1}+M_{X_{i+1}})$-Mori fiber space$/Z$ $\phi_{i+1}: X_{i+1}\rightarrow S_{i+1}$,
        \item two $\Rr$-Cartier $\Rr$-divisors $L_{i+1}$ and $H_{i+1}$ on $X_{i+1}$, 
            \item two real number $0\leq l_{i+1}\leq l_i$ and $h_i\leq h_{i+1}\leq 1$, and
            \item a Sarkisov link$/Z$ $f_i: X_i\dashrightarrow X_{i+1}$ as in \textbf{Case 1}, or \textbf{Case 2.1}, or \textbf{Case 2.2}, or \textbf{Case 3.1}, or \textbf{Case 3.2}, or \textbf{Case 3.3} of Theorem \ref{thm: scaling sarkisov},
    \end{itemize}
    such that
     \begin{itemize}
        \item $(W,B_W+l_{i+1}L_W+h_{i+1}H_W+M_W)\geq (X_{i+1},B_{i+1}+l_{i+1}L_{i+1}+h_{i+1}H_{i+1}+M_{X_{i+1}})$,
        \item $B_{i+1},L_{i+1}$ and $H_{i+1}$ are the birational transforms of $B_{i},L_{i}$ and $H_{i}$ on $X_{i+1}$ respectively,
        \item $K_{X_{i+1}}+B_{i+1}+l_{i+1}L_{i+1}+h_{i+1}H_{i+1}+M_{X_{i+1}}\sim_{\Rr,S_{i+1}}0$, and
        \item $K_{X_{i+1}}+B_{i+1}+l_{i+1}L_{i+1}+h_{i+1}H_{i+1}+M_{X_{i+1}}$ is nef$/Z$.
    \end{itemize}
    Notice that the assumptions hold when $i=0$.
\item[\textbf{Step 3}] If $l_{i+1}=0$, we stop and let $n:=i+1$. Otherwise, we replace $i$ with $i+1$ and return to \textbf{Step 2}.
\end{itemize}
The following diagram gives the birational maps and Mori fiber spaces in this construction:
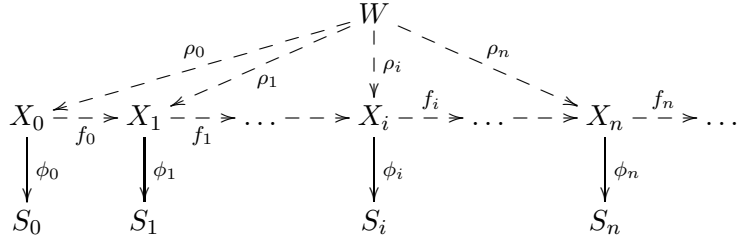
\begin{center}$\xymatrix{
 & & &W\ar@{-->}[dlll]_{\rho_0}\ar@{-->}[dll]^{\rho_1}\ar@{-->}[d]^{\rho_i}\ar@{-->}[drr]^{\rho_n}& & & \\
      X_0 \ar@{-->}[r]_{f_0}\ar@{->}[d]^{\phi_0}&X_1\ar@{-->}[r]_{f_1}\ar@{->}[d]^{\phi_1}&\dots\ar@{-->}[r]&X_i\ar@{-->}[r]^{f_i}\ar@{->}[d]^{\phi_i}&\dots\ar@{-->}[r]&X_n\ar@{->}[d]^{\phi_n}\ar@{-->}[r]^{f_n}& \dots \\
    S_0 & S_1 & & S_i& &S_n & }$
\end{center}
\end{cons}

\begin{lem}
Assumptions and notation as in Construction \ref{cons: sarkisov scaling}. Then
\begin{enumerate}
    \item there are only finitely many possibilities of $X_i$ up to isomorphism, and
    \item the Sarkisov program of $(X,B+M_X)$ with scaling of $(L_W,H_W)$ terminates, i.e. there exists an integer $n>0$ such that $l_n=0$.
\end{enumerate}
\end{lem}
\begin{proof}
We construct two a ample$/Z$ $\Rr$-divisors $\Gamma_W$ and $N_W$ on $W$, boundary divisors $\Delta_{i,W}$ for every $i\gg 0$, and a finite dimensional affine subspace $\mathcal{V}$ of $\Weil_{\Rr}(W)$ in the following way:
\begin{itemize}
    \item Assume that $h_k>0$ for some integer $k\geq 0$. In this case, since $H_W$ is big and nef$/Z$, we may write $H_W\sim_{\Rr,Z}C_W+G_W$, such that $C_W$ is ample$/Z$, $G_W$ is effective, and $(W,B_W+L_W+H_W+C_W+G_W+M_W)$ is g-terminal. We let $\Gamma_W\sim_{\Rr,Z}M_W+\frac{h_k}{2}C_W$ be a general ample$/Z$ $\Rr$-divisor on $W$ such that 
    $$(W,B_W+l_iL_W+\frac{h_k}{2}C_W+h_iG_W+\Gamma_W+(h_i-h_k)C_W)$$ 
    is klt for every $i\geq k$. We define $N_W:=\frac{h_k}{2}C_W$, $\Delta_{i,W}:=B_W+l_iL_W+\frac{h_k}{2}C_W+h_iG_W+\Gamma_W+(h_i-h_k)C_W$ for every $i\geq 0$, and $\mathcal{V}$ be the finite dimensional affine subspace of $\Weil_{\Rr}(W)$ generated by the irreducible components of $B_W,L_W,G_W,\Gamma_W$ and $C_W$.
    \item Assume that $h_i=0$ for every integer $i\geq 0$. In this case, $l_i=1$ for every $i\geq 0$. Since $L_W$ is big and nef$/Z$, we may write $L_W\sim_{\Rr,Z}D_W+J_W$, such that $D_W$ is ample$/Z$, $J_W$ is effective, and $(W,B_W+L_W+D_W+J_W+M_W)$ is g-terminal.  We let $\Gamma_W\sim_{\Rr,Z}M_W+\frac{1}{2}D_W$ be a general ample$/Z$ $\Rr$-divisor on $W$ such that 
     $$(W,B_W+\frac{1}{2}D_W+J_W+\Gamma_W)$$ 
    is klt for every $i\geq k$. We define $N_W:=\frac{1}{2}D_W$, $\Delta_{i,W}:=B_W+\frac{1}{2}D_W+J_W+\Gamma_W)$ for every $i\geq 0$, and $\mathcal{V}$ be the finite dimensional affine subspace of $\Weil_{\Rr}(W)$ generated by the irreducible components of $B_W,J_W$ and $\Gamma_W$.
\end{itemize}
Since  $K_{X_i}+B_i+l_iL_i+h_iH_i+M_{X_i}$ is nef$/Z$ for every $i$ and $K_{W}+B_W+l_iL_W+h_iH_W+M_W\sim_{\Rr,Z}K_W+\Delta_{i,W}$, $\rho_i: W\dashrightarrow X_i$ is a weak log canonical model$/Z$ of $(W,\Delta_{i,W})$. Since $\Delta_{i,W}\in\mathcal{L}_{N_W}(\mathcal{V})$ for every $i$, we deduce (1) by applying Theorem \ref{thm: finiteness ltm}.

Suppose that (2) does not hold. Then $l_i>0$ for every $i>0$. By (1), there exists a strictly increasing sequence $\{i_j\}_{j=1}^{+\infty}$, such that $X_{i_j}\cong X_{i_k}$ for any $j,k\geq 1$.  By construction, we have $l_{i_j}=l_{i_k}$ and $h_{i_j}=h_{i,k}$ for every $j,k\geq 1$. Thus $l_i=l_{i_1}$ and $h_i=h_{i_1}$ for every $i\geq i_1$. In particular, $f_i$ is a Sarkisov link$/Z$ of as in \textbf{Case 2} or \textbf{Case 3} of Theorem \ref{thm: scaling sarkisov} for every $i\geq i_1$.

Suppose that $f_i$ is a Sarkisov link$/Z$ as in \textbf{Case 3} of Theorem \ref{thm: scaling sarkisov} for some $i\geq i_1$. By Lemma \ref{lem: sar scaling terminalization}(3), $f_{i+1}$ is a Sarkisov link$/Z$ as in \textbf{Case 3} of Theorem \ref{thm: scaling sarkisov}, which implies that $f_j$ is a Sarkisov link$/Z$ as in \textbf{Case 3} of Theorem \ref{thm: scaling sarkisov} for every $j\geq i$. By Theorem \ref{thm: scaling sarkisov}(6) and Lemma \ref{lem: sar scaling terminalization}(4), $f_j$ is not a Sarkisov link$/Z$ as in \textbf{Case 3.1} of Theorem \ref{thm: scaling sarkisov} for any $j\geq i$. Since $\rho(X_j)>0$ for every $j$, by Lemma \ref{lem: sar scaling terminalization}(5)(6.a), for every $j\gg i$,  $f_j$ is a Sarkisov link$/Z$ as in \textbf{Case 3.3} of Theorem \ref{thm: scaling sarkisov}. This contradicts to Lemma \ref{lem: sar scaling terminalization}(6).

Therefore, we may assume that $f_j$ is a Sarkisov link$/Z$ as in \textbf{Case 2} of Theorem \ref{thm: scaling sarkisov} for every $j\geq i_1$. Since $\rho(X_j)\leq\rho(W)$ for every $j$, by Lemma \ref{lem: sar scaling terminalization}(1)(2.a), for every $j\gg i_1$, $f_j$ is a Sarkisov link$/Z$ as in \textbf{Case 2.2} of Theorem \ref{thm: scaling sarkisov}. This contradicts to Lemma \ref{lem: sar scaling terminalization}(2).
\end{proof}

\section{Proof of the main theorem}
\begin{thm}\label{thm: existence sarkisov precise}
Assume that 
\begin{itemize}
    \item $W\rightarrow Z$ is a contraction between normal varieties,
    \item $(W,B_W+M_W)$ is a $\mathbb Q$-factorial gklt g-pair with associated nef$/Z$ $\bb$-divisor $M$, such that $K_W+B_W+M_W$ is not pseudo-effective$/Z$,
    \item $\rho_X: W\dashrightarrow X$ and $\rho_Y: W\dashrightarrow Y$ are two $(K_W+B_W+M_W)$-non-positive maps$/Z$ such that $(\rho_X)_*(K_W+B_W+M_W)=K_X+B_X+M_X$ and $(\rho_Y)_*(K_W+B_W+M_W)=K_X+B_Y+M_Y$,
    \item $\phi_X: X\rightarrow S_X$ is a $(K_X+B_X+M_X)$-Mori fiber space$/Z$ and $\phi_Y: X\rightarrow S_Y$ is a $(K_Y+B_Y+M_Y)$-Mori fiber space$/Z$.
\end{itemize}
\begin{center}$\xymatrix{
 & W\ar@{-->}[dl]_{\rho_X}\ar@{-->}[dr]^{\rho_Y}& \\
      X \ar@{->}[d]_{\phi_X}\ar@{-->}[rr]^{f}   &  & Y\ar@{->}[d]^{\phi_Y} \\
    S_X & &S_Y }$
\end{center}
Then
\begin{enumerate}
    \item the induced birational map $f: X\dashrightarrow Y$ is given by a finite sequence of Sarkisov links$/Z$, i.e. $f$ can be written as $X_0\dashrightarrow X_1\dots\dashrightarrow X_n\cong Y$, where each $X_{i}\dashrightarrow X_{i+1}$ is a Sarkisov link,
    \item $f$ is a Sarkisov program$/Z$ with double scaling, i.e., each $X_{i}\dashrightarrow X_{i+1}$ is a Sarkisov link$/Z$ as in Construction \ref{cons: sarkisov scaling},
    \item for any real number $\epsilon>0$, if $(W,B_W+M_W)$ is generalized $\epsilon$-lc, then $(X_i,B_i+M_{X_i})$ is generalized $\epsilon$-lc for every $i$, where each $B_i$ is the birational transform of $B_W$ on $X_i$.
\end{enumerate} 
\end{thm}
\begin{proof}
By Lemma \ref{lem: terminalization and mmp}, possibly replacing $(W,B_W+M_W)$ with a common log resolution of $(W,B_W+M_W)$, $(X,B_X+M_X)$ and $(Y,B_Y+M_Y)$, we may assume that $(W,B_W+M_W)$ is g-terminal, $\rho_X$ and $\rho_Y$ are morphisms, and $M=M_W$. By our assumptions, $-(K_X+B_X+M_X)$ is $\phi_X$-ample and $-(K_Y+B_Y+M_Y)$ is $\phi_Y$-ample. Therefore, there are
\begin{itemize}
    \item two ample$/Z$ $\Rr$-divisors $A_{S_X}$ and $A_{S_Y}$ on $S_X$ and $S_Y$ respectively, 
    \item two general ample$/Z$ $\Rr$-divisors $L_X$ and $H_Y$ on $X$ and $Y$ respectively, and
    \item two general big and nef$/Z$ $\Rr$-divisors $L_W$ and $H_W$ on $W$,
\end{itemize}
such that
\begin{itemize}
    \item $L_X\sim_{\Rr,Z}-(K_X+B_X+M_X)+\phi_X^*A_{S_X}$,
    \item $H_Y\sim_{\Rr,Z}-(K_Y+B_Y+M_Y)+\phi_X^*A_{S_Y}$, and
    \item $(X,B_X+L_X+M_X)$ and $(Y,B_Y+H_Y+M_Y)$ are both gklt,
    \item $L_W=\rho_X^*L_X=(\rho_X^{-1})_*L_X$ and $H_W=\rho_X^*L_X=(\rho_Y^{-1})_*L_X$, and
    \item $(W,B_W+2(H_W+L_W)+M_W)$ is g-terminal.
\end{itemize}
By Lemma \ref{lem: sar scaling terminalization}, we may let $f': X:=X_0\dashrightarrow X_1\dots\dashrightarrow X_n$ be a Sarkisov program of $(X,B_X+M_X)$ with scaling of $(L_W,H_W)$ as in Construction \ref{cons: sarkisov scaling}. 

We show that $f=f'$ and $h_n=1$. Indeed, if $h_n<1$, then since $(\rho_Y)_*(K_W+B_W+h_nH_W+M_W)\sim_{\Rr,S_Y}-(1-h_n)H_Y$, $\phi_Y$ is a $(K_Y+B_Y+h_nH_Y+M_Y)$-Mori fiber space$/Z$, which implies that $K_W+B_W+h_nH_W+M_W$ is not pseudo-effective$/Z$. However, since $K_{X_n}+B_{n}+h_nH_{n}+M_{X_n}$ is nef$/Z$ by construction, $K_W+B_W+h_nH_W+M_W$ is pseudo-effective$/Z$, a contradiction.

Let $D_n$ be a general ample$/Z$ $\Rr$-divisor on $X_n$ and $D_Y,D_W$ the birational transforms of $D_n$ on $Y$ and $W$ respectively. Since $D_Y$ is big$/Z$, $D_Y$ is ample$/S_Y$. Pick $0<\epsilon\ll 1$ and let $\Delta_W\sim_{\Rr,Z}B_W+H_W+M_W+\epsilon D_n$ be an effective $\Rr$-divisor on $W$ such that $(W,\Delta_W)$ is klt. Then $W\dashrightarrow X_n$ and $W\dashrightarrow Y$ are both log canonical models$/Z$ of $(W,\Delta_W)$, which implies that $X_n\cong Y$.

Let $\Sigma_n$ be a $\phi_n$-vertical curve, then 
$$\phi_Y^*A_{S_Y}\cdot\Sigma_n=(K_{Y}+B_{Y}+H_{Y}+M_{Y})\cdot\Sigma_n=0,$$
which implies that $\Sigma_n$ is $\phi_Y$-vertical. Thus $\phi_n$ and $\phi_Y$ define the same Mori fiber space and the theorem follows.
\end{proof}

\begin{proof}[Proof of Theorem \ref{thm: existence generalized Sarkisov link}]
Since a $(K_W+B_W+M_W)$-MMP$/Z$ is a $(K_W+B_W+M_W)$-non-positive map$/Z$, Theorem \ref{thm: existence generalized Sarkisov link} follows from Theorem \ref{thm: existence sarkisov precise}.
\end{proof}

\end{document}